\documentclass[a4paper]{paper}
\usepackage{multirow}
\usepackage{mathptmx}
\usepackage{amsmath,amsthm,amsfonts,mathtools}
\usepackage{graphicx}
\usepackage{subcaption}
\usepackage{derivative}
\usepackage{algorithm}
\usepackage{algpseudocode}
\DeclareMathOperator{\diag}{diag}

\usepackage{listings}
\usepackage{subdepth}

\newcommand{\mi}{\mathrm{i}}
\newcommand{\be}{\begin{equation}}
\newcommand{\ee}{\end{equation}}
\newcommand{\ba}{\begin{array}}
\newcommand{\ea}{\end{array}}

\newcommand{\comment}[1]{}

\usepackage{pgfplots, pgfplotstable, tikz}
\pgfplotstableset{col sep=comma}
\pgfplotsset{compat=1.16}
\usepackage{hyperref}
\newtheorem{theorem}{Theorem}

\title{Numerical Solution of Nonclassical Boundary Value Problems}

\author{Boito, Paola \thanks{The author is a member of GNCS-INDAM.  The author acknowledges the MIUR Excellence Department Project awarded to the Department of Mathematics, University of Pisa, CUP I57G22000700001.}\\
  Dipartimento di Matematica\\
  Universit\`a di Pisa \\
  \texttt{paola.boito@unipi.it}
  \and 
  Eidelman, Yuli \\
  School of Mathematical Sciences\\
  Raymond and Beverly Sackler Faculty of
Exact Sciences\\
Tel-Aviv University\\
\texttt{eideyu@tauex.tau.ac.il}
\and 
  Gemignani, Luca  \thanks{The author is a member of GNCS-INDAM. The author is  partially supported 
 by  European Union - NextGenerationEU under the National Recovery and Resilience Plan (PNRR) - Mission 4 Education and research - Component 2 From research to business - Investment 1.1 Notice Prin 2022 - DD N. 104  2/2/2022, titled Low-rank Structures and Numerical Methods in Matrix and Tensor Computations and their Application, proposal code 20227PCCKZ – CUP I53D23002280006  and by the Spoke 1 ``FutureHPC \& BigData”  of the Italian Research Center on High-Performance Computing, Big Data and Quantum Computing (ICSC)  funded by MUR Missione 4 Componente 2 Investimento 1.4: Potenziamento strutture di ricerca e creazione di "campioni nazionali di R\&S (M4C2-19 )" - Next Generation EU (NGEU).}\\
  Dipartimento di Informatica\\
  Universit\`a di Pisa \\
  \texttt{luca.gemignani@unipi.it}
}

\begin{document}
\maketitle

\begin{abstract}
We provide a new approach to obtain  solutions of  certain evolution equations set in a Banach space and
equipped with nonlocal boundary conditions.    From this approach we derive a family of numerical schemes for the  approximation of the solutions. We show by numerical tests  that these schemes are numerically  robust and computationally efficient.
\end{abstract}

\section{Introduction}
In this work we focus on the solution of a class of linear differential problems of the form 
$$\frac{dv}{dt}=Av,\quad 0<t<T,$$
set in a Banach space, with the nonlocal integral condition
$$\frac1{T}\int_0^Tv(t)\;dt=f,$$
where the function $f$ and the linear, closed, possibly unbounded operator $A$ are given. 

For instance, $A$ could be chosen as a second derivative w.r.t.~to a space variable $x$, yielding the familiar form of the homogeneous heat equation, albeit with less common boundary conditions. The study of the heat equation with integral boundary conditions goes back to Cannon \cite{cannon1963solution}; the existence and properties of the solution are investigated in \cite{tikhonov2003uniqueness, popov2005exponential}, see also \cite{kalenyuk2012problem} and references therein.  

In a finite-dimensional setting, $A$ is a linear operator acting on a finite-dimensional space, that is, a matrix. The solution to our differential problem can then be expressed as the action of a function of $A$ on $f$, namely, $v(t)=\psi_t(A)f$, where $\psi_t(z)=\frac{Tze^{zt}}{e^{zT}-1}$. Note in passing that this function is closely related to the reciprocal of the $\varphi_1$ function, which is of interest in the development of exponential integrators: see e.g., \cite{BEG2, botchev2021residual, HO, jimenez2020efficient} and references therein. 

A discussion of the finite-dimensional problem was proposed in \cite{BEG1} and \cite{BEG2}. More precisely, \cite{BEG1} introduced an expression for $v(t)$ based on a partial fraction decomposition of the function $\psi_t(z)$, thus motivating the development of a structured algorithm for the solution of shifted quasiseparable linear systems. It also presented an acceleration technique for this rational decomposition, which introduced a cubic polynomial term. In \cite{BEG2}, generalization of this idea led to a family of mixed polynomial-rational approximations of  $v(t)$, where the polynomial part is in fact a Bernoulli polynomial of arbitrary degree.  In addition to its theoretical interest, this approximation formula allows for the design of effective numerical methods for the computation of $v(t)$, particularly when $A$ is a structured matrix. Let us also mention that the approach based on Bernoulli polynomials for the solution of differential problems with nonlocal conditions was suggested for the first time in \cite{eidelman2017application}. 

For the more general case where the differential problem is set in a Banach space, we prove the existence and uniqueness of the solution $v(t)$ and characterize it
via a family of  mixed polynomial-rational expansions w.r.t.~the operator $A$. Each expansion contains a purely polynomial term of arbitrary degree, which is related to the Bernoulli polynomials, followed by a series of rational terms. The evaluation of each rational term involves the inversion of $A$ plus a multiple of the identity operator. 

From this result we derive a general numerical procedure for computing an approximation of $v(t)$ up to a given tolerance (Algorithm 1). An interesting feature of this approach is the fact that successive rational terms can be computed independently: this allows us to fine-tune the accuracy of the approximation by adding further terms as needed, without the need to recompute the whole approximation. Moreover, in order to improve the efficiency of the implementation, one may also employ parallelization and/or structured methods for the inversion of families of shifted operators; a finite-dimensional example in presence of quasiseparable structure is found in \cite{BEG1}.   

Numerical tests highlight the effectiveness of this approach and suggest strategies for a good choice of the degree of the polynomial term and of the number of rational terms. After providing results for the matrix case, we focus on a model problem of parabolic equation. For this problem we investigate the behavior of the mixed approximation in combination with two different approaches: the classical method of lines, based on a finite-difference semi-discretization in space, and a ``functional'' approach where the expansion is applied directly to the infinite-dimensional operator. Each rational term is then computed as the solution of a boundary value problem. Among the proposed tests, this functional approach is the most innovative and promising application of our mixed approximation formula. While more computationally demanding, it circumvents the numerical obstacles encountered in the finite-dimensional discretization and proves to be more robust and flexible.  

The paper is organized as follows. Section \ref{sec:maintheorem} recalls the problem under study and contains the main theoretical result of the paper, namely, the statement and proof for the mixed polynomial-rational expansion of the solution (Theorem \ref{main}). This is followed by Algorithm 1, which computes an approximation of $v(t)$. Sections \ref{sectwo} and \ref{secthree} are devoted to numerical tests for the matrix case and for the parabolic model equation, respectively. Section \ref{sec:conclusions} summarizes our contributions and sets out ideas for future work, particularly  for further investigation and improvement of the functional approach.

\section{The Abstract Nonlocal Differential Problem}\label{sec:maintheorem}
\setcounter{equation}{0}

Let $X$ and $Y$ be two Banach spaces.  Consider the  linear homogeneous differential equation
\be\label{l1}
\frac{dv}{dt}=Av,\quad 0<t<T,
\end{equation}
with the boundary conditions
\be\label{l2}
\frac1{T}\int_0^Tv(t)\;dt=f.
\end{equation}
Here $A \colon D(A)\rightarrow X$ is a linear unbounded closed operator with the domain 
$D(A)\subset X$,  and $f\in D(A)$. By solution of the problem (\ref{l1}), (\ref{l2}) we mean a continuous function $v \colon [0,T]\rightarrow X$
on $[0,T]$ with values in $X$ and with $v(t)\in D(A),\;0<t<T$,  such that 
$v\in \mathcal C^1((0,T);X)$ and $Av\in \mathcal C((0,T);X)$, and (\ref{l1}), (\ref{l2}) hold.
Without loss of generality one can assume that $T =2\pi$. 
Recall that 
a linear  operator $T\colon X\rightarrow Y$ 
 is closed if whenever $x_k\rightarrow x$
 in $X$ 
 and $Tx_k\rightarrow y$
 in $Y$, we have $T x=y$. We say that $\lambda \in \mathbb C$ is a regular point of $A$ if $A-\lambda I$, $I=I_X$,  is bijective  and $(A-\lambda I_X)^{-1}$ is a bounded operator, i.e., 
 \[
\parallel (A-\lambda I)^{-1} \parallel = \sup\{\frac{\parallel (A-\lambda I)^{-1} y \parallel}{\parallel y\parallel}\colon y\neq 0\}<+\infty.
\]
The following result establishes the basis for the design of numerical schemes for  computing a numerical approximation of  the solution $v(t)$.  A preliminary version of this result  in the finite dimensional case  has first appeared  in \cite{BEG1,BEG2}.

\begin{theorem}\label{main}
Assume that all the complex numbers
$$
\mu_k=\mi k,\quad k=\pm1,\pm2,\dots
$$
are regular points of the operator $A$ and there is a constant $C>0$ such that
\begin{equation}\label{gul}
\|(A-\mu_kI)^{-1}\|\le\frac{C}{|k|},\quad k=\pm1,\pm2,\dots.
\end{equation} 
Also, suppose that $f\in D(A^2)$.
Then the problem (\ref{l1}), (\ref{l2}) has a unique solution which is given
by the formula
\begin{equation}\label{base}
v(t)=f+(t-\pi) Af + 2\left(\sum_{k=1}^{\infty}\Sigma_k(Af)\cos kt +
\sum_{k=1}^{\infty}\Upsilon_k(Af)\sin kt\right)
\end{equation}
with
\begin{equation}\label{shiftsys}
\left\{ \begin{array}{ll}
\Sigma_k=A(A+\mi k  I )^{-1}(A-\mi  k I )^{-1},\\
\Upsilon_k=k^{-1} A^2 (A+\mi k  I )^{-1}(A-\mi k  I )^{-1}, 
\end{array}\right. \quad k=1,2,\dots
\end{equation}
Moreover for any $n\ge0$ under the additional assumption $f\in D(A^{2n+2})$ the
 solution $v(t)$ satisfies the formula
\begin{equation}\label{ll}
    v(t)=p_n(t)+s_n(t)
\end{equation}
with 
\begin{equation}\label{llk}
p_n(t)=\sum_{k=0}^{2n+1} \frac{(2\pi)^k}{k!} B_k(\frac{t}{2\pi})A^kf,
\end{equation}
and
\begin{equation}\label{nate}
s_n(t)=(-1)^{n}2\sum_{k=1}^{\infty}\frac1{k^{2n}}
(\Sigma_k(A^{2n+1}f)\cos kt+ \Upsilon_k(A^{2n+1}f)\sin kt), 
\end{equation}
where $B_m(t)$ are the well-known Bernoulli polynomials:
$$
B_m(t)=\sum_{j=0}^m\frac1{j+1}\sum_{k=0}^j(-1)^k
\left(\begin{array}{c}j\\k\end{array}\right)(t+k)^m.
$$
\end{theorem}

\begin{proof}

For any $f\in D(A^2)$ we consider the sequence of partial sums
\be\label{gal}
v_n(t)=f+(t-\pi) Af + 2\left(\sum_{k=1}^{n}\Sigma_k(Af)\cos kt +
\sum_{k=1}^{n}\Upsilon_k(Af)\sin kt\right),\;n=0,1,2,\dots.
\end{equation}
At first we prove that the sequence $v_n(t)$ in (\ref{gal}) 
converges uniformly in $t$ on $[0,2\pi]$. Indeed using the first formula in 
(\ref{shiftsys}) we have
$$
\|\Sigma_k(Af)\cos(kt)\|\le 
\|(A+\mi k  I )^{-1}\|\|(A-\mi  k I )^{-1}\|\|A^2f\|,\;k=0,1,2,\dots\;
0\le t\le 2\pi
$$
and from (\ref{gul}) follows the uniform convergence in $t\in[0,2\pi]$ for the first
series in (\ref{gal}). Using the second formula in (\ref{shiftsys}) we
get
$$
\Upsilon_k(Af)\sin kt=k^{-1} A(A+\mi k  I )^{-1}(A-\mi k  I )^{-1}\sin kt
(A^2f),\;k=1,2,\dots\;0\le t\le 2\pi.
$$
Using the equality $k^{-1}A=k^{-1}(A-\mi k  I )+\mi I$ we get
\begin{gather*}
\Upsilon_k(Af)\sin kt=k^{-1} (A+\mi k  I )^{-1}(A^2f)\sin kt+
\mi(A^2+k^2I)^{-1}(A^2f)\sin kt,\;0\le t\le 2\pi,\\
k=1,2,\dots\;
\end{gather*}
which implies the uniform convergence in $t\in[0,2\pi]$ for the second
series in (\ref{gal}). Thus there is a continuous function $v(t)$ on $[0,2\pi]$
which is a uniform limit of the sequence $v_n(t)$ on $[0,2\pi]$. Moreover using 
the equalities
$$
\int_0^{2\pi}\cos kt\;dt=\int_0^{2\pi}\sin kt\;dt=0,\;k=1,2,\dots,\quad
\int_0^{2\pi}(t-\pi)\;dt=0 
$$
we get
$$
\int_0^{2\pi}v_n(t)\;dt=2\pi f,\;n=0,1,2,\dots
$$
and passing to the limit as $n\to\infty$ we obtain the condition (\ref{l2}).

Now consider the sequences
\begin{equation}\label{gal1}
Av_n(t)=Af+(t-\pi) A^2f + 2\left(\sum_{k=1}^{n}\Sigma_k(A^2f)\cos kt +
\sum_{k=1}^{n}\Upsilon_k(A^2f)\sin kt\right),\;n=0,1,2,\dots
\end{equation}
and
\begin{equation}\label{gal3}
\frac{dv_n}{dt}=Af + 2\left(\sum_{k=1}^{n}\Sigma_k(Af)(-k)\sin kt +
\sum_{k=1}^{n}\Upsilon_k(Af)k\cos kt\right),\;n=0,1,2,\dots
\end{equation}
Comparing (\ref{gal1}) with (\ref{gal3}) and using (\ref{shiftsys}) we get
$$
\frac{dv_n}{dt}+\phi(t)A^2f=Av_n(t)+\phi_n(t)A^2f
$$
with
$$
\phi(t)=t-\pi,\quad \phi_n(t)=-\sum_{k=1}\frac2{k}\sin kt
$$
One can easily see that $\phi_n(t),\;n=1,2,\dots$ is the sequence of the 
partial sums of the Fourier series of the function $\phi(t)$ on the
segment $[0,2\pi]$. Hence it follows that $\phi_n(t)$ converges to $\phi(t)$ as
$n\to\infty$ uniformly in $t$ on $[\delta,\gamma]$
for any $\delta,\gamma$ with $0<\delta<\gamma<2\pi$.
One should prove that the sequences $Av_n(t)$ in (\ref{gal1}) and
 $\frac{dv_n}{dt}$ in (\ref{gal3}) converge uniformly in $t$ on 
$[\delta,\gamma]$. To this end one can check the uniform  convergence in $t$ on
$[\delta,\gamma]$ of the sums 
\begin{equation}\label{ok}
\sum_{k=1}^{n}\Sigma_k\cos kt,\;\sum_{k=1}^{n}\Upsilon_k\sin kt.
\end{equation}
We apply the Abel transformation to these sums. We have
$$
\sum_{k=1}^n\Sigma_k\cos kt=
\sum_{k=1}^{n-1}(\Sigma_k-\Sigma_{k+1})\sum_{j=1}^k\cos jt+
\Sigma_n\sum_{j=1}^n\cos jt.
$$
Using the formula
$$
\sum_{j=1}^k \cos jt=-\frac1{2}+\frac1{2}\frac{\sin\frac{2k+1}{2}t}
{\sin \frac{t}{2}}
$$
we get
$$
\sum_{k=1}^n\Sigma_k \cos kt=-\frac1{2}\sum_{k=1}^{n-1}(\Sigma_k-\Sigma_{k+1}) + \frac1{2}\sum_{k=1}^{n-1}(\Sigma_k-\Sigma_{k+1})
\frac{\sin\frac{2k+1}{2}t}{\sin \frac{t}{2}}-\frac1{2}\Sigma_n+\frac1{2}\Sigma_n
\frac{\sin\frac{2n+1}{2}t}{\sin \frac{t}{2}}
$$
which implies
$$
\sum_{k=1}^n\Sigma_k \cos kt=-\frac1{2}\Sigma_1+\frac1{2}\Sigma_n
\frac{\sin\frac{2n+1}{2}t}{\sin\frac{t}{2}}+
\frac1{2}\sum_{k=1}^{n-1}(\Sigma_k-\Sigma_{k+1})
\frac{\sin\frac{2k+1}{2}t}{\sin \frac{t}{2}}.
$$

Similarly for  the second series we have
$$
\sum_{k=1}^n\Upsilon_k\sin kt=
\sum_{k=1}^{n-1}(\Upsilon_k-\Upsilon_{k+1})\sum_{j=1}^k\sin jt+
\Upsilon_n\sum_{j=1}^n\sin jt.
$$
Using the formula
$$
\sum_{j=1}^k \sin jt=\frac1{2}\cot\frac{t}{2}-
\frac1{2}\frac{\cos\frac{2k+1}{2}t}
{\sin \frac{t}{2}}
$$
we get
$$
\sum_{k=1}^n\Upsilon_k \sin kt=
\frac1{2}\cot\frac{t}{2}\sum_{k=1}^{n-1}(\Upsilon_k-\Upsilon_{k+1})-
\frac1{2}\sum_{k=1}^{n-1}(\Upsilon_k-\Upsilon_{k+1})
\frac{\cos\frac{2k+1}{2}t}
{\sin \frac{t}{2}}+
$$
$$
\frac1{2}\cot\frac{t}{2}\Upsilon_n-
\frac1{2}\Upsilon_n\frac{\cos\frac{2n+1}{2}t}{\sin\frac{t}{2}}
$$
which implies
$$
\sum_{k=1}^n\Upsilon_k \sin kt=\frac1{2}\cot\frac{t}{2}\Upsilon_1-
\frac1{2}\Upsilon_n
\frac{\cos\frac{2n+1}{2}t}{\sin\frac{t}{2}}-
\frac1{2}\sum_{k=1}^{n-1}(\Upsilon_k-\Upsilon_{k+1})
\frac{\cos\frac{2k+1}{2}t}{\sin \frac{t}{2}}.
$$
In order to prove convergence of \eqref{ok}, we need to study the behavior of $\Sigma_k$, $\Upsilon_k$, $\Sigma_k-\Sigma_{k+1}$ and $\Upsilon_k-\Upsilon_{k+1}$ as $k\rightarrow\infty$.
Using (\ref{shiftsys}) we have
\begin{eqnarray*}
&& \Sigma_k-\Sigma_{k+1}=A(A^2+k^2I)^{-1}-A(A^2+(k+1)^2I)^{-1}\\
&& = A(A^2+k^2I)^{-1}(A^2+(k+1)^2I)^{-1}[(A^2+(k+1)^2I)-(A^2+k^2I)]\\
&&=A(2k+1)(A^2+k^2I)^{-1}(A^2+(k+1)^2I)^{-1}
\end{eqnarray*}
and analogously
$$
\Upsilon_k-\Upsilon_{k+1}=
\frac{1}{k(k+1)}(A^2+(3k^2+3k+1)I)(A^2+k^2I)^{-1}(A^2+(k+1)^2I)^{-1}.
$$
Now observe that from the identity
$$
\frac{A}{k}= \frac1k(A-ikI)+i\,I
$$
it follows 
\begin{equation}\label{auxeq}
    \frac{A}{k}(A-ikI)^{-1}=\frac1k+i(A-ikI)^{-1}.
\end{equation}
Using the inequalities (\ref{gul}) together with \eqref{auxeq} we get
\begin{equation}\label{gl}
\|\Sigma_k-\Sigma_{k+1}\|,
\|\Upsilon_k-\Upsilon_{k+1}\|\le \frac{K}{k^2},\quad k=1,2,\dots,
\end{equation}
where $K$ is a constant.
Moreover using (\ref{shiftsys}) we have
$$
\Sigma_k=(A+ikI)^{-1}+\mi k(A^2+k^2I)^{-1},\;
\Upsilon_k=\frac1{k}-k(A^2+k^2I)^{-1}.
$$
This by virtue of (\ref{gul}) implies
\begin{equation}\label{ug}
\lim_{n\to\infty}\Sigma _n=\lim_{n\to\infty}\Upsilon_n=0.
\end{equation}
Finally for any $\delta,\gamma$ with $0<\delta<\gamma<2\pi$ we have
\begin{equation}\label{bsin}
\left|\sin\frac{t}{2}\right|\ge c>0,\;0<\delta\le t\le\gamma<2\pi.
\end{equation}
Thus combining the relations (\ref{gl}), (\ref{ug}), (\ref{bsin}) together we 
conclude that the sums (\ref{ok}) converge uniformly in $t$ on 
$[\delta,\gamma]$. Hence the same is true for the sequences (\ref{gal1}), (\ref{gal3}).
Since $A$ is a closed operator,  we obtain that $Av_n(t)\rightarrow A v(t)$.  Analogously, the uniform convergence of the derivatives implies that the derivative of the limit is the limit of the derivatives, and, hence, $\displaystyle\frac{dv_n}{dt} \rightarrow \displaystyle\frac{dv}{dt}$. Therefore, we conclude that 
$v(t)$ is a solution of the equation (\ref{l1}).

Proceeding in the same way as in \cite{BEG1,BEG2} with the assumption 
$f\in D(A^{2n+2})$ we obtain the formulas 
(\ref{ll}), (\ref{llk}), (\ref{nate}) for the solution $v(t)$.

Let us now prove the uniqueness of the solution. Let $v(t)$ be a solution of the 
 linear homogeneous differential problem
\begin{equation}\label{l1bis}
\frac{dv}{dt}=Av,\; 0<t<T,
\end{equation}
\begin{equation}\label{hh}
\frac1{2\pi}\int_0^{2\pi}v(t)\;dt=0.
\end{equation}
Note that by integrating \eqref{l1bis} from $0$ to $2\pi$ and applying \eqref{hh} we obtain
\begin{equation}\label{v02pi}
    v(2\pi)-v(0)=\int_0^{2\pi}\frac{dv}{dt}dt=A\int_0^{2\pi}v(t) dt=0.
\end{equation}
The function $v(t)$ is continuous on $[0,2\pi]$ and has Fourier coefficients
$c_k,\;k=0,\pm1,\pm2,\dots$ given by
$$
c_k=\frac1{2\pi}\int_0^{2\pi}v(t)e^{-\mu_kt}\;dt,\quad k=\pm1,\pm2,\pm
$$
and
\be\label{el}
c_0=\frac1{2\pi}\int_0^{2\pi}v(t)\;dt=0.
\end{equation}
Multiplying the equation (\ref{l1bis}) by $\frac1{2\pi}e^{-\mu_kt}$ and 
integrating from $0$ to $2\pi$ we get
$$
\frac1{2\pi}\int_0^{2\pi}\frac{dv}{dt} e^{-\mu_kt}\;dt=
A\frac1{2\pi}\int_0^{2\pi}v(t)e^{-\mu_kt}\;dt,\quad k=\pm1,\pm2,\dots.
$$
Integrating by parts in the left hand side and applying \eqref{v02pi} we get
$$
\mu_kc_k=Ac_k,\;k=\pm1,\pm2,\dots.
$$
which implies
\be\label{ele}
c_k=0,\;k=\pm1,\pm2,\dots
\end{equation}
Thus from (\ref{ele}), (\ref{el}) it follows that all the Fourier coefficients
of the continuous function $v(t)$ are zeros. Hence $v(t)=0,\;0\le t\le 2\pi$.
\end{proof}

Based on this theorem, we get a  family of  polynomial/rational approximations  of the solution of \eqref{l1}, \eqref{l2} of the form:
\begin{equation}\label{sch1}
v(t)\simeq v_{n,\ell}(t)=p_n(t)+s_{n,\ell}(t)
\end{equation}
where $s_{n,\ell}(t)$ is a partial sum of the series (\ref{nate}), i.e.
$$
s_{n, \ell}(t)=(-1)^{n}2\sum_{k=1}^{\ell}\frac1{k^{2n}}
(\Sigma_k(A^{2n+1}f)\cos kt+ \Upsilon_k(A^{2n+1}f)\sin kt)
$$
and 
\[
v(t)=v_{n,\ell}(t) + r_{n, \ell}(t)
\]
with residual 
\[
r_{n, \ell}(t)=(-1)^{n}2\sum_{k=\ell+1}^{+\infty}\frac1{k^{2n}}
(\Sigma_k(A^{2n+1}f)\cos kt+ \Upsilon_k(A^{2n+1}f)\sin kt).
\]
The bulk of the paper deals with theoretical and computational issues associated with the computation of $v_{n,\ell}(t)$.   The results of a preliminary numerical experience  reported in \cite{BEG1}, \cite{BEG2}  and \cite{LG1} upon the finite dimensional case  promote
some considerations concerning the selection of the parameters $n$ and $\ell$. 
\begin{enumerate}
\item In the finite dimensional case  it is shown that  the larger the degree of $p_n(t)$ is,  the better the convergence of the series of $v(t)$  is.  However,  in finite precision arithmetic large values  of $n$  generally lead to stability  issues in the computation of the approximant \eqref{sch1}. Moreover,  in view of Theorem  \ref{main} large values of $n$ imply additional smoothness requirements  upon the function $f(t)$ in \eqref{l2}. Therefore,  in our numerical experience  the better strategy is to set the value of $n$ as small as possible, typically $n\in \{0,1,2,3,4\}$, and then determine the value of $\ell$ so as to obtain the desired accuracy. 
\item An adaptive  technique for the selection of the value of $\ell$  in the finite dimensional case is presented in \cite{LG1}.  For a fixed  $t\in [0,T]$ once the value of $n$  and a threshold tolerance $\epsilon>0$ are given,  then we can  determine the value of $\ell=\ell(t)$ by imposing the condition 
\[
\parallel v_{n, \ell-1}(t)-v_{n, \ell}(t)\parallel/\parallel v_{n, \ell}(t)\parallel \leq \epsilon.
\]
In this way  the approximation $v_{n, \ell}(t)$  is constructed incrementally
by  adding one term  at a time until a  fixed tolerance is reached. 
\end{enumerate}
According to these facts, we consider
the following general procedure for computing an approximation $v_{n, \ell}(t)$  of $v(t)$ for $t\in [0,T]$. 

\begin{algorithm}\label{a1}
\caption{}
    \begin{algorithmic}[1] 
        \State \textbf{Select} the values of $n \in \{0,1,2,3,4\}$ and $tol$. 
        \State \textbf{Define} a coarse  grid of points $\mathcal I_S=\{t_0, \ldots, t_m\}\subset[0,T]$.  
        \State \textbf{Compute},  for any $t_i\in \mathcal I_S$, $v_{n, \ell_i}(t_i)$ as in \eqref{sch1} where $\ell_i$ is the minimum 
index such that 
\begin{equation}\label{errtest}
\parallel v_{n, \ell_i-1}(t_i)-v_{n, \ell_i}(t_i)\parallel_\infty/\parallel v_{n, \ell_i}(t_i)\parallel_\infty \leq tol
\end{equation}
holds.  
\State  \textbf{Set} $\widehat \ell=\max_{1\leq i\leq m} \ell_i$  be the length of our rational expansion.
            \State \textbf{return} $v_{n,\widehat \ell}(t)$  as the approximation of $v(t)$ over $[0,T]$.
    \end{algorithmic}
\end{algorithm}

The computational effort of  \textbf{Algorithm 1} at step 3  basically  amounts to 
determine $\Sigma_k(Af)$ and $\Upsilon_k(Af)$ for $k=1, \ldots, \widehat \ell$.
The accuracy of the approximation $v_{n,\widehat \ell}(t)$ returned as output by 
\textbf{Algorithm 1}  is measured by introducing a finer grid of points 
$\mathcal I_L\subset [0,T]$ and then setting 
\begin{equation}\label{errmeas}
err=\max_{t_i\in \mathcal I_L} \parallel v(t_i)-v_{n, \widehat \ell}(t_i)\parallel_\infty/\parallel v(t_i)\parallel_\infty .
\end{equation}

 The selection of grid points can be problem dependent  or related to specific properties of the solution function. For testing purposes we   generally make use of equispaced grid of points.
Theoretical estimates of the norm of the residual $r_{n,\ell}(t)$ in terms of of the quantities  \eqref{errtest} are currently unavailable. The feasibility and robustness of  \textbf{Algorithm 1} is validated by numeric simulations shown in the next sections.

\section{The matrix case}\label{sectwo}
\setcounter{equation}{0}

Here it is assumed that $X={\mathbb C}^N$, $A$ is an $N\times N$ matrix, $f$ is
an $N$-dimensional column and $v(t)={\rm col}(v_i(t))_{i=1}^N$ is a vector
function. 
For testing purposes  let us suppose that $A$ is diagonalizable by a unitary congruence,  that is, $A= Q^H D Q$ with $Q^H Q=I_N$ and $D=\diag[\lambda_1, \ldots, \lambda_N]$. Then it can be proved  by direct calculations that the solution of \eqref{l1},\eqref{l2}  can be expressed as 
\begin{equation}\label{practeq}
v(t)=Q^H \diag[e^{\lambda_1t} \psi_1(\lambda_1 2 \pi), \ldots, e^{\lambda_N t} \psi_1(\lambda_N 2 \pi)] Q  f,
\end{equation}
where $\psi_1(z)=z/(e^z-1)$  is a meromorphic function with poles $z_k=2 \pi \mu_k= 2 \pi \mi k,\quad k=\pm1,\pm2,\dots$. It is tacitly assumed that the computation of $v(t)$ by means of \eqref{practeq} gives the "exact" solution of \eqref{l1},\eqref{l2}.  Thus, such solution   can be  used in our numerical simulations to provide error estimates according to \eqref{errmeas}.

We have implemented the computation of the approximation  
 $v_{n, \ell}(t)$  of  $v(t)$    as described in \eqref{sch1}  
 and \textbf{Algorithm 1}.  The coarse and fine grid of points are formed by equispaced 
nodes in the interval $[0,T]$. The input parameters are $m, n\in \mathbb N$ and $tol\in \mathbb R$.  Set $\mathcal I_S=\{t_i=2\pi(i-1)/(m-1)\colon 1\leq i\leq m\}$ the  set of $m$  uniformly spaced points in the interval $[0, T]=[0, 2 \pi]$.  The finer grid of points is determined so that  $\mathcal I_L=\{t_i=2\pi(i-1)/(m-1)^2\colon 1\leq i\leq (m-1)^2+1\}$  includes $\mathcal I_s$.

Our test suite consists of the following set of  rank-structured matrices:
\begin{enumerate}
 \item{\tt  Random unitary  Hessenberg matrices.} Test problems
with known eigenvalues are generated as follows. A unitary diagonal matrix $D$ is generated and its eigenvalues noted. A unitary matrix $Q$, random with respect to Haar measure, is generated, and the random unitary matrix 
$B = Q^H D Q$  formed. Then $B$  is  transformed to upper Hessenberg form by unitary congruence,   that is,  
$P^H B P=H $,  to yield an upper Hessenberg
unitary matrix  $H$ with known eigenvalues, which is then factored into the 
form $H= (QP)^H D  (QP)$.  The matrix $V=(QP)^H$ is the  unitary eigenvector  matrix. 
\item {\tt Banded circulant matrices.}  We consider  banded circulant matrices of the form 
 \[
 C=\left[\begin{array}{ccccc}  a_0 &&&a_2 &a_1\\
 a_1& \ddots &&&a_2\\a_2& \ddots & \ddots \\
 & \ddots &\ddots &\ddots \\
 &&a_2 & a_1 & a_0 \end{array}\right] \in \mathbb C^{N\times N}. 
 \]
 We recall that  $C$ can be diagonalized by a  Discrete Fourier Transform.  Specifically,  if $F=(f_{j,\ell})$, $f_{j, \ell}=\displaystyle\frac{e^{\mi 2 \pi (j-1)(\ell-1)/N}}{\sqrt{N}}$, $1\leq j, \ell \leq N$, is the Fourier matrix of order $N$ then  $F^H C F=\diag[\lambda_1, \ldots, \lambda_N]$ with 
 \[
 \lambda_k=\sum_{j=0}^2 a_j e^{-\mi 2 \pi (k-1) j/N}, \quad 1\leq k\leq N.
 \]
\item {\tt Shifted and scaled 1D Laplacian matrices.}  We consider  $N\times N$ tridiagonal Toeplitz matrices  generated by 
\begin{equation}\label{ss1D}
A=\sigma   \ {\tt gallery('tridiag'}, N, 1, -2, 1) -\gamma \ {\tt eye(N)}
\end{equation}
with $\sigma, \gamma \in \mathbb R^+$.  These matrices can be diagonalized  by means of a  discrete sine transform  $E=(\sqrt{2/(N+1)} \sin \displaystyle \frac{\pi i j}{N+1})_{i,j=1}^N$.
\end{enumerate}
 The figures and tables below show  the results of our numerical experiments. 

\begin{figure}[ht] 
\begin{subfigure}[b]{0.5\linewidth}
    \centering
    \includegraphics[width=0.75\linewidth]{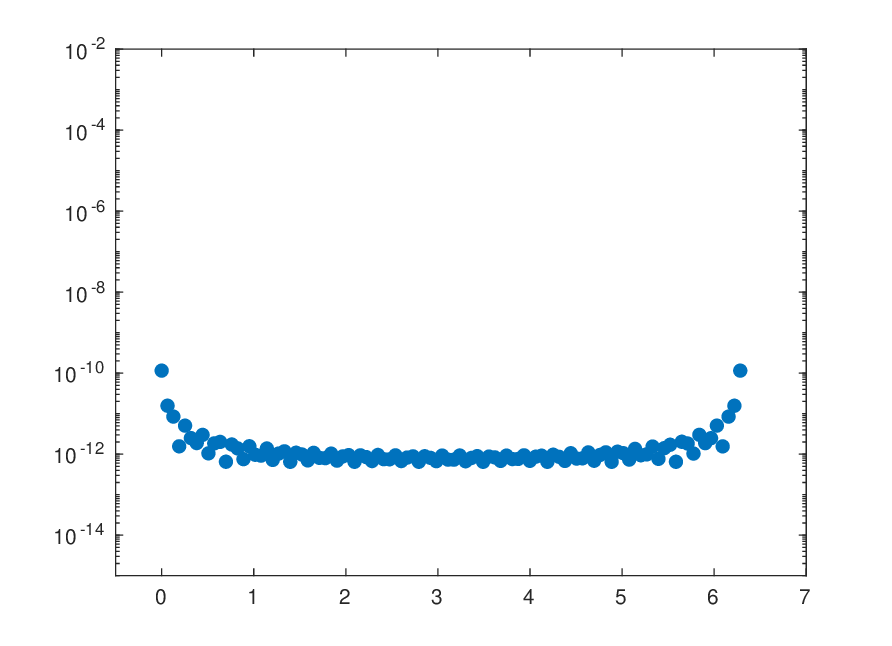} 
    \caption{Plot of computed error for $n=1$.} 
    \label{5_1} 
    \vspace{4ex}
  \end{subfigure}
  \begin{subfigure}[b]{0.5\linewidth}
    \centering
    \includegraphics[width=0.75\linewidth]{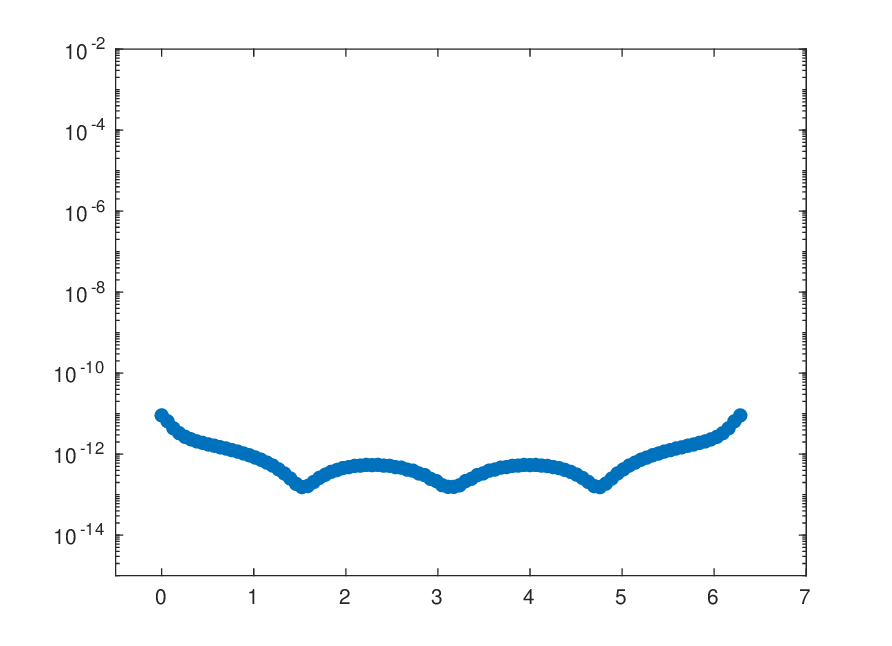} 
    \caption{Plot of computed error for $n=2$.}  
    \label{5_2} 
    \vspace{4ex}
  \end{subfigure} 
  \begin{subfigure}[b]{0.5\linewidth}
    \centering
    \includegraphics[width=0.75\linewidth]{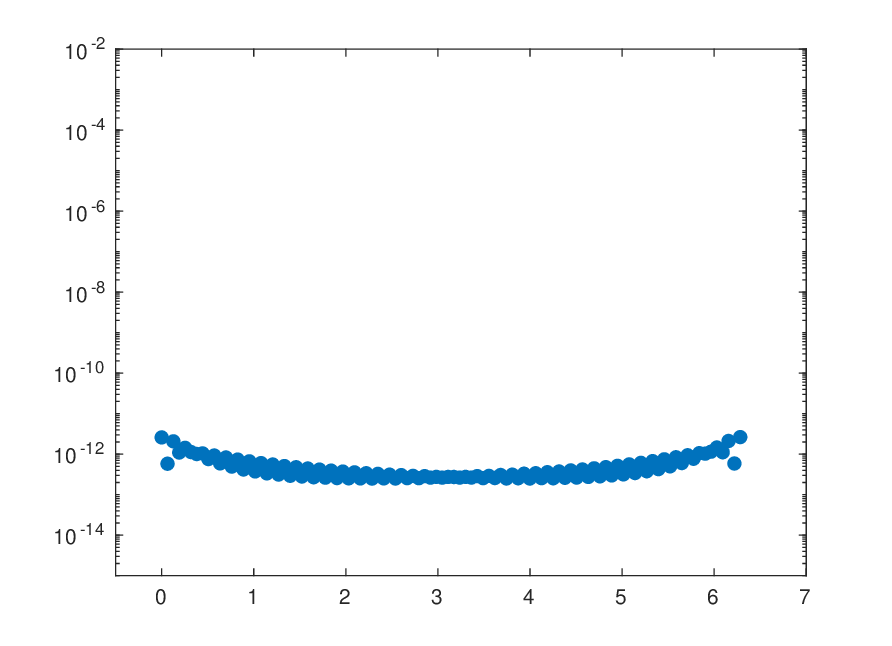} 
      \caption{ Plot of computed error for $n=3$.} 
    \label{5_3} 
  \end{subfigure}
  \begin{subfigure}[b]{0.5\linewidth}
    \centering
    \includegraphics[width=0.75\linewidth]{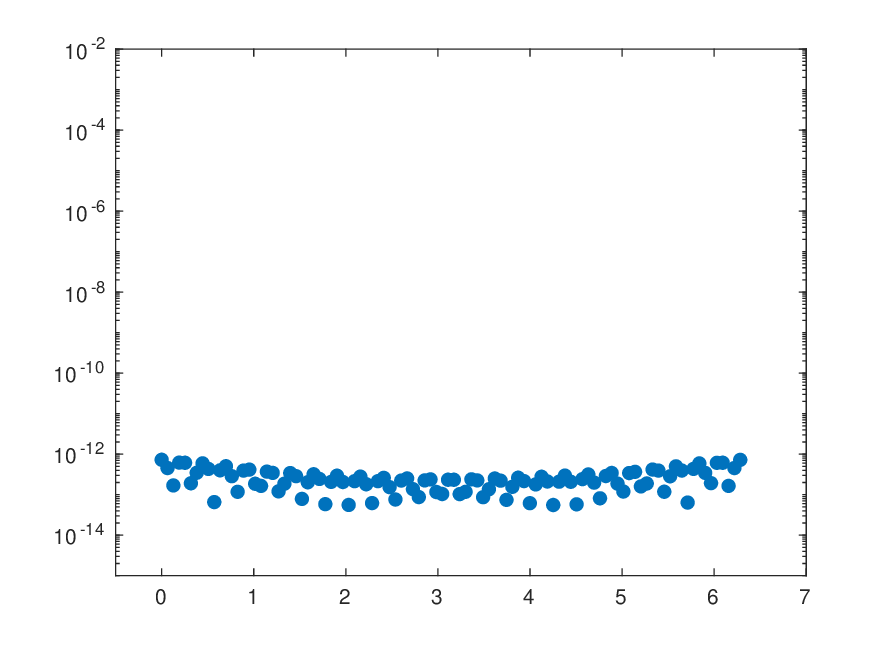} 
    \caption{Plot of computed error for $n=4$.} 
    \label{5_4} 
  \end{subfigure} 
  \caption{Illustration of   the  measured error  \eqref{errmeas} for Test 1  with $tol=1.0e\!-12$ and $m=10$  for different values of $n\in \{1,2,3,4\}$. }
  \label{fig5} 
\end{figure}

Random unitary Hessenberg matrices provide an easy test class.  The eigenvalues are located on the unit circle in the complex plane,  the norm of $A$ is bounded by 1 and the solution of  the  shifted linear systems associated with the computation of $\Sigma_k(Af)$ and $\Upsilon_k(Af)$ in \eqref{base}, \eqref{shiftsys}
 is   generally well conditioned.  Therefore,   our proposed  schemes  perform quite well in this case.  
In Figure \ref{fig5} we illustrate the plots of the measured error over the grid $\mathcal I_L$ for $N=1024$, $tol=1.0e\!-12$ and different values of $n\in \{1,2,3,4\}$. Notice that the error is flattened out as $n$ increases and, moreover, the stopping tolerance gives  a quite precise  measure of the norm of the residual. We guess that the error  behaviour is determined by  the decreasing of the algorithmic error in the computation of the  mixed polynomial/rational approximation of the differential problem.

\begin{figure}
  \centering
  \subfloat[]{\includegraphics[width=0.4\textwidth]{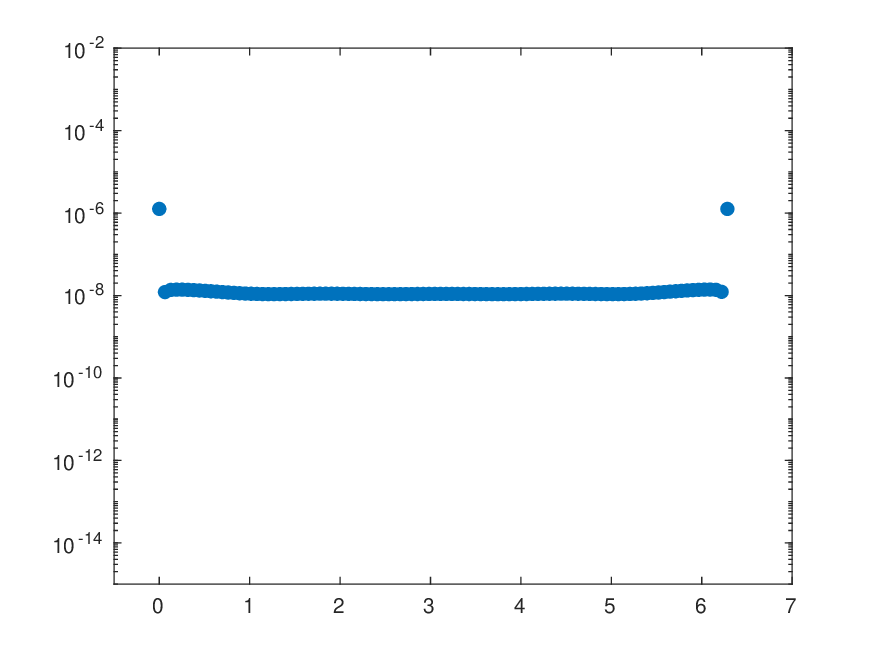}}
  \hfill
  \subfloat[]{\includegraphics[width=0.4\textwidth]{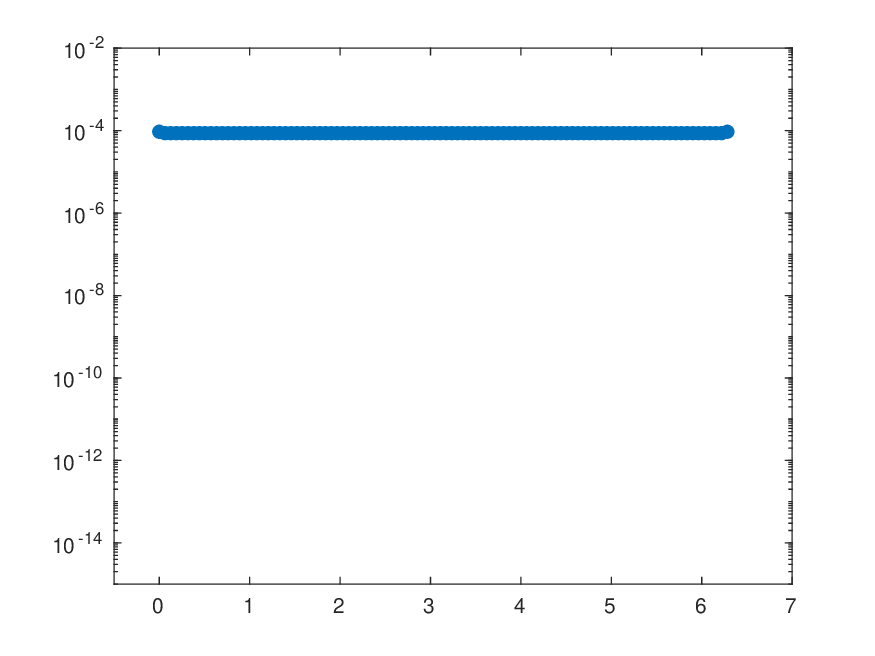}}
  \caption{Illustration of   the  measured error  \eqref{errmeas} for Test 2  with $tol=1.0e\!-12$, $m=10$, $n=1$, $a_0=a_1=1, a_2=1+1.0e\!-4$ in (a) and $a_0=a_1=1, a_2=1+1.0e\!-6$ in (b). }
\label{fig30}
\end{figure}

\begin{figure}[ht] 
  \begin{subfigure}[b]{0.5\linewidth}
    \centering
    \includegraphics[width=0.75\linewidth]{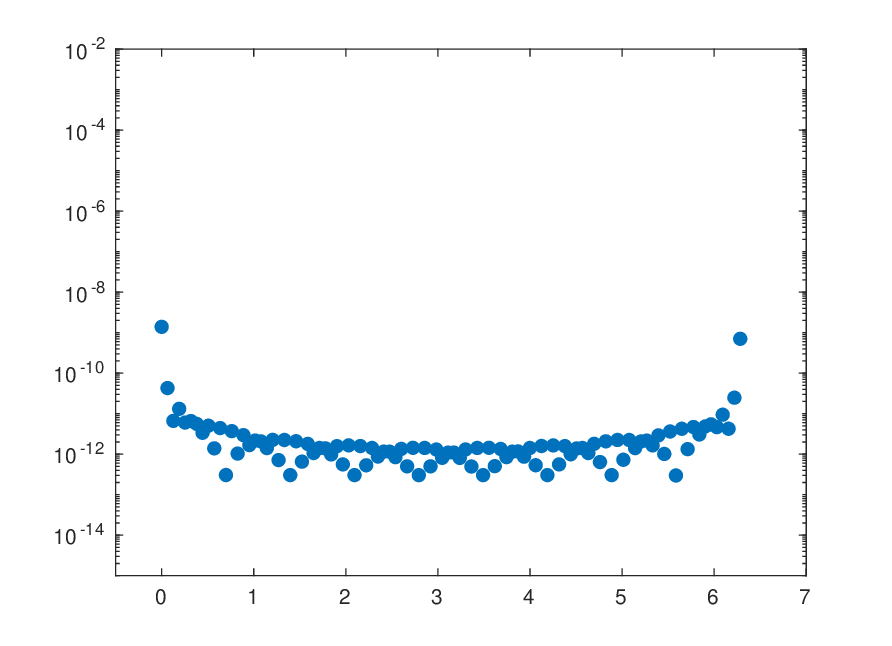} 
    \caption{Plot of computed error for $n=1$.} 
    \vspace{4ex}
  \end{subfigure}
  \begin{subfigure}[b]{0.5\linewidth}
    \centering
    \includegraphics[width=0.75\linewidth]{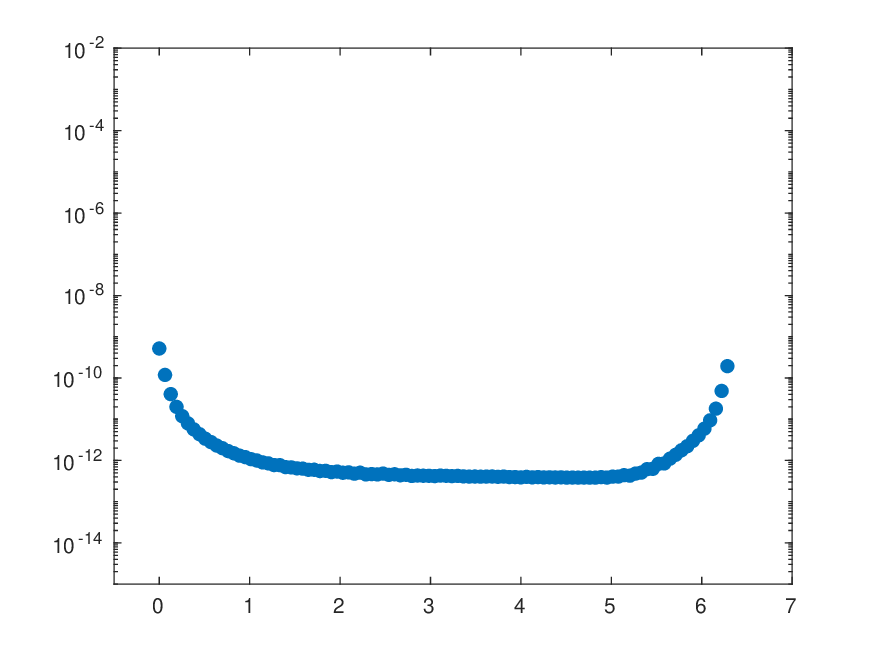} 
    \caption{Plot of computed error for $n=2$.}  
    \vspace{4ex}
  \end{subfigure} 
  \begin{subfigure}[b]{0.5\linewidth}
    \centering
    \includegraphics[width=0.75\linewidth]{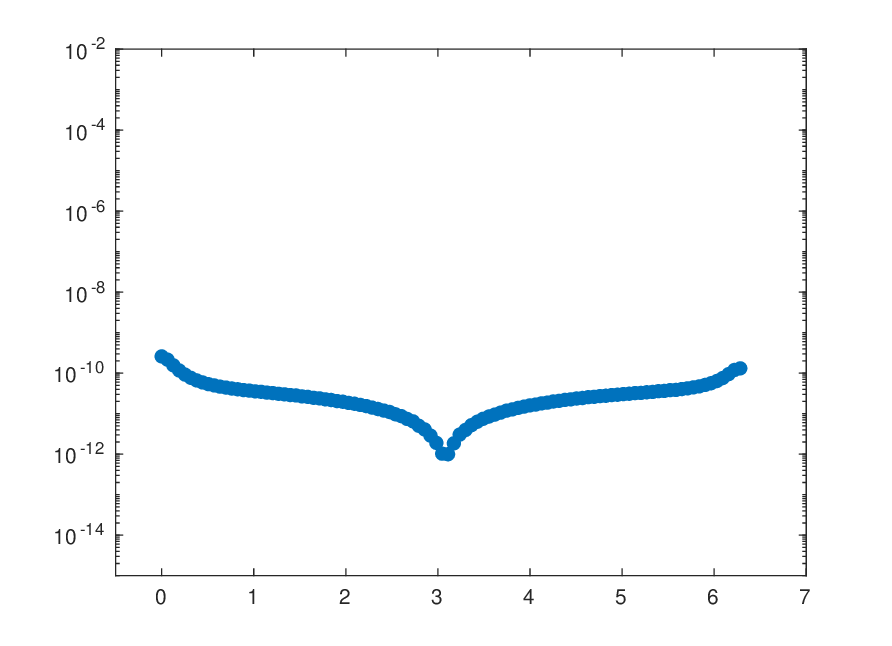} 
      \caption{ Plot of computed error for $n=3$.} 
  \end{subfigure}
  \begin{subfigure}[b]{0.5\linewidth}
    \centering
    \includegraphics[width=0.75\linewidth]{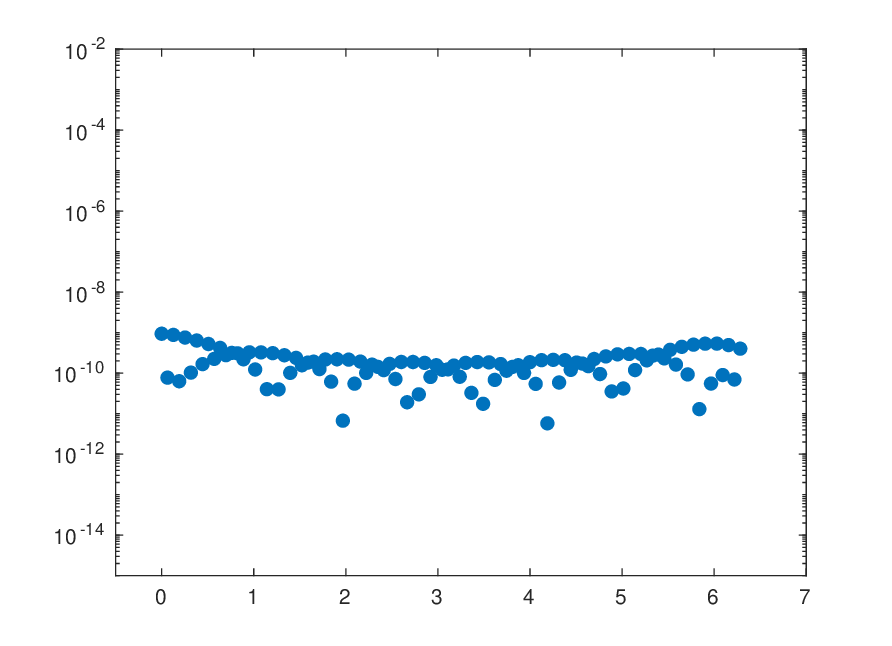} 
    \caption{Plot of computed error for $n=4$.} 
  \end{subfigure} 
  \caption{Illustration of   the  measured error  \eqref{errmeas} for Test 1 with $a_0=a_1=a_2=1$, $tol=1.0e\!-12$ and $m=10$  for different values of $n\in \{1,2,3,4\}$. }
  \label{fig1new} 
\end{figure}

To experience with more difficult tests we consider banded circulant matrices.  For 
$a_0=a_1=a_2=1$ and $N=2k$ the matrix  $C$ has  one eigenvalue equals to $\mi$.   In view of Theorem \ref{main} the solution  of \eqref{l1},\eqref{l2} does not exist or is not unique.  The MATLAB implementation of our code reports "ill-conditioning" warnings  in Matlab's  backslash command  used to solve the shifted linear systems involved in the computation of the rational approximant $s_n(t)$. 
Perturbed values of the coefficients make it possible to tune the conditioning of the shifted linear systems. 
In  Figure \ref{fig30} we show the plots of the computed error for $N=1024$, $tol=1.0e\!-12$, $n=1$ and $a_0=a_1=1, a_2=1+1.0e\!-4$ and 
$a_0=a_1=1, a_2=1+1.0e\!-6$, respectively.  This figure clearly indicates the impact of the conditioning of the shifted linear systems 
on the overall accuracy of the computed approximation. 
Differently, for $N=2k-1$  the proposed scheme performs quite satisfactorily. In Figure \ref{fig1new} we show  the measured error   for Test 1 with $a_0=a_1=a_2=1$, $tol=1.0e\!-12$ and $m=10$.   Also,  the acceleration of the convergence due to the increasing value of $n$ is dramatic. For $n=0$  the stopping criterion  \eqref{errtest} is not fulfilled within 50000 iterations.   Differently,  for $n=4$ we find $\widehat \ell=22$. 

\begin{figure}
  \centering
  \includegraphics[width=0.6\textwidth]{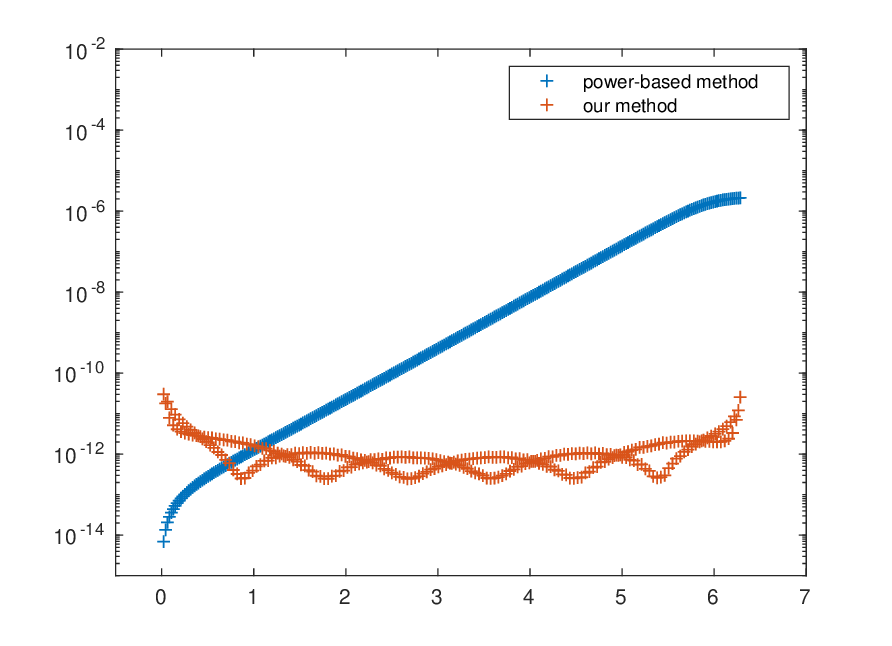}
  \caption{ Error plots  generated by  the power-based method and our algorithm applied or solving Test 1 with $a_0=a_1=a_2=1$, $tol=1.0e\!-12$, $n=2$  and $m=30$.}
\label{ff200}
\end{figure}

It is worth mentioning that \eqref{practeq}   can be rewritten as  $v(t)={\tt expm}(A t) v(0)$ where ${\tt expm}(A)$ is the matrix exponential of $A$. This relation provides an easy way based on the power method to give approximations of the solution vectors $v(t_i)$, $t_i\in \mathcal I_L$. Assume  that $B={\tt expm}(A t_1)$ and  $v(0)$
are available. Then we can  determine  the remaining vectors $v(t_i)$ by setting 
$v(t_{i+1})=B v(t_i)$, $i\geq 0$. In Figure \ref{ff200} we compare the accuracy of this power-based method with our algorithm applied for solving Test 1 with $a_0=a_1=a_2=1$, $tol=1.0e\!-12$, $n=2$  and $m=30$.  The matrix $B$ is computed by means of the MATLAB function {\tt expm}.

\begin{figure}
  \centering
  \subfloat[]{\includegraphics[width=0.4\textwidth]{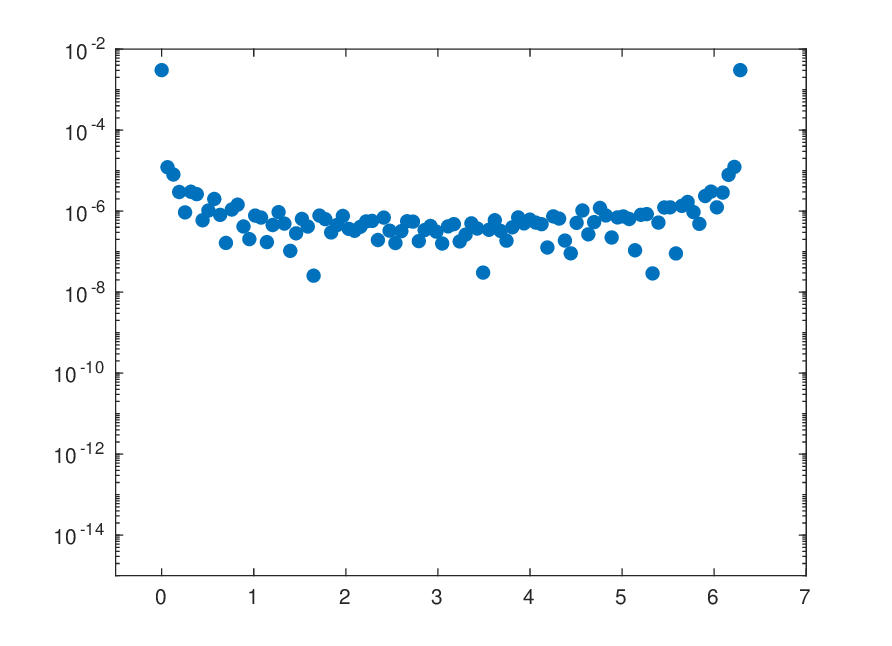}}
  \hfill
  \subfloat[]{\includegraphics[width=0.4\textwidth]{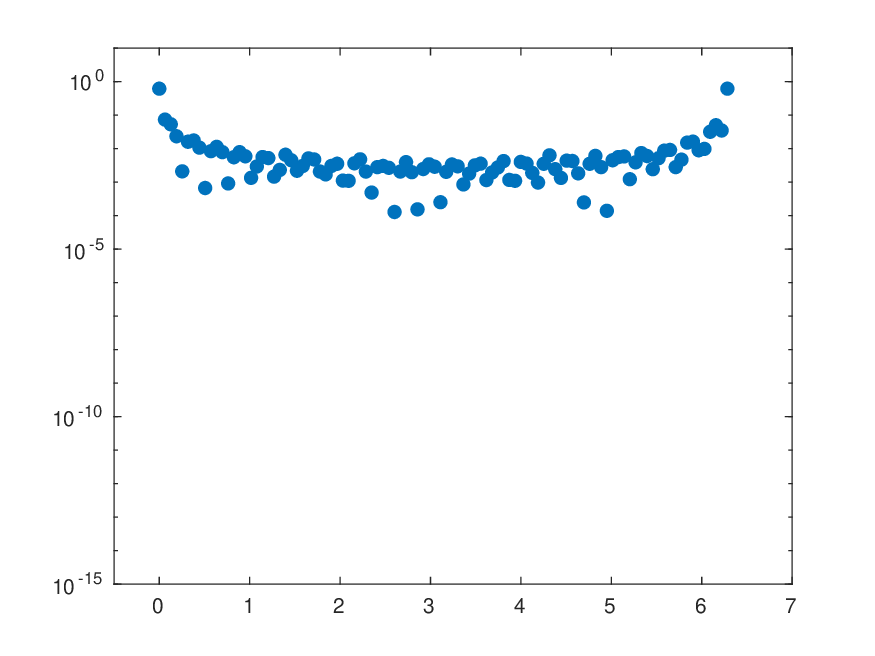}}
  \caption{Illustration of   the  measured error  \eqref{errmeas} for Test 3  with $\sigma=1$, $\gamma=100$,  $tol=1.0e\!-12$, $m=10$, $n=1$ in (a) and $n=2$ in (b). }
\label{fig40}
\end{figure}

Matrices $A\in \mathbb R^{N\times N}$ of the form \eqref{ss1D} 
are useful to check the dependence of the convergence rate and the accuracy  w.r.t.~the distribution of the eigenvalues  and the norm of $A$.  
This  set has already been considered for testing purposes in \cite{BEG1,BEG2}.  For $\sigma=1$ and $\gamma=0$  all the proposed schemes with $n\in \{0,1,2,3,4\}$ perform well in terms of accuracy  confirming the results in \cite{BEG1}.  For $N=1024$ and $tol=1.0e\!-12$
the scheme with $n=0$ is unpractical since it requires more than 50000 terms to satisfy the stopping criterion. For $n=1, 4$ we find $\widehat l=1714$ and $\widehat l=21$, respectively, by showing again the  efficiency of the  convergence  acceleration introduced by the Bernoulli polynomials.  Increasing the value of $\sigma$  or $\gamma$ introduces some numerical difficulties.   For $\sigma=1$ and $\gamma=100$ the shifted linear systems are  still well-conditioned but  both the magnitude of eigenvalues and the norm of $A$ increase.  As a side effect we observe a deterioration of the convergence and the accuracy of  our algorithm. For $n=1, 2$ we find $\widehat l=8896$ and $\widehat l=457$, respectively. In  Figure \ref{fig40} we show the plots of the computed error for  shifted 1D Laplacian matrices  with $N=1024$, $tol=1.0e\!-12$,   and $n=1, 2$. The schemes with $n=3,4$ return inaccurate results.  Similar behaviors for the error are found for 
$\sigma=100$ and $\gamma=0$.

In conclusion,  numerical experiments with the finite dimensional case reveal that the computation of   the approximant $v_{n, \ell}(t)$ given 
in \eqref{sch1}  can be prone to numerical inaccuracies essentially due to: 
\begin{enumerate}
\item a possibly large norm of the matrix  $A$  which makes  unreliable the use of acceleration schemes based on Bernoulli polynomials; 
\item the effects of a potential ill-conditioning of the shifted linear systems involved in the computation of the trigonometric expansion.
\end{enumerate}
 Moreover, the occurrence of large eigenvalues in magnitude   determines two additional  effects: 
\begin{enumerate}
\item the convergence generally deteriorates as more and more  terms are needed into the expansion to evaluate the function at the eigenvalues; 
\item  since the convergence is slow the stopping criterion  can be   satisfied even the approximation is far from the  "exact" value. 
\end{enumerate}
The analysis of these weak points  motivates the  extension of the results in  \cite{BEG1,BEG2} to the functional setting given in Theorem \ref{main}.  In particular, in the next section we will show that  the extension  makes possible the application of our approach for solving  certain differential problems with nonlocal boundary conditions. 

\section{The model one-dimensional problem for parabolic  equation}\label{secthree}
\setcounter{equation}{0}

Consider the  differential problem
\be\label{al}
\left\{\begin{array}{ll}
u_t=\sigma u_{xx}+c(x)u,&0<x<1,\;0<t< T,\\
u(0,t)=u(1,t)=0,&0<t< T.\end{array}
\right.
\end{equation}
with  the nonlocal condition 
\be\label{lor}
\frac1{T}\int_0^Tu(x,t)\;dt=f(x),\quad f\in D(A).
\end{equation}

In principle, several  numerical schemes for computing  the solution of \eqref{al},  \eqref{lor} can be devised based on  Theorem \ref{main}. 
Specifically, we treat the problem (\ref{al}) as a problem for an ordinary  
differential equation in the Banach space $X=L^2(0,1)$ of the form 
\begin{equation}\label{l11} 
\frac{dv}{dt}=Av,\quad 0<t<T, 
\end{equation} 
Here $A=\sigma \displaystyle\frac{d^2}{dx^2}+c(x)$, where the domain $D(A)$ consists of the functions $w(x)$ with first derivative  $w'$ absolutely continuous on $[0,1]$ such that $w''\in L^2(0,1)$, and satisfying the boundary conditions $w(0)=w(1)=0$. It is worth noting that expansion \eqref{base} requires $f\in D(A^2)$, where $D(A^2)$ is the set of functions $w(x)$ on $[0,1]$ such that $w'''(x)$ is absolutely continuous, $w^{(4)}(x)\in L^2(0,1)$, and $w(0)=w(1)=w''(0)=w''(1)=0$.

Notice that if we take $c(x)=0$ for simplicity, then for any  
$f\in D(A^k)$ the values $A^kf$ are in fact the even-order derivatives of the given 
function $f$, namely $A^kf=f^{(2k)}$. Similarly, if $c(x)$ is constant, then $A^kf$ can be expressed as linear combinations of the even-order derivatives of the given function $f(x)$. 
So we can assume that the elements 
$g_j=A^jf,\;j=0,\dots,k$ are given as well. Using Theorem \ref{main}, with the 
 assumption $f\in D(A^{2n+3})$ , we have that the solution $v(t)$ satisfies the formula 
\begin{equation}\label{l1l} 
v(t)=v_{n, \ell}(t) +r_{n,\ell}(t)=p_n(t)+s_{n,l}(t)  + r_{n,\ell}(t)
\end{equation} 
with 
\begin{equation}\label{l1k} 
p_n(t)=\sum_{k=0}^{2n+1} \frac{(2\pi)^k}{k!} B_k\left(\frac{t}{2\pi}\right)g_k, 
\end{equation} 
and 
\be\label{nater} 
s_{n,l}(t)=(-1)^{n}2\sum_{k=1}^{l}\frac1{k^{2n}} 
\left(\Sigma_k A^2g_{2n-1}\cos kt+ \Upsilon_k A^2g_{2n-1}\sin kt\right). 
\end{equation} 
We transform the formula (\ref{nater}) to an equivalent form. We define the operators 
\be\label{irav} 
V_k=(A^2+k^2I)^{-1},\quad k=1,2,\dots . 
\end{equation} 
Using (\ref{shiftsys}) we have  
that (\ref{nater}) takes the form \be\label{nater1} 
s_{n,l}(t)=(-1)^{n}2\sum_{k=1}^{l}\frac1{k^{2n}}\left( 
V_kg_{2n+2}\cos kt+\frac1{k}V_kg_{2n+3}\sin kt\right). 
\end{equation} 
So the core of the algorithm is the evaluation of the quantities 
$$ 
V_kg=(A^2+k^2I)^{-1}g,\; k=1,2,\dots,\quad g\in X. 
$$ 
The computation of these values may be done as follows. 
\begin{enumerate} 
\item Compute the solution $p^{(k)}$ of the equation  
\be\label{uri1} 
(A-(ki)I)p^{(k)}=g. 
\end{equation} 
\item Compute the solution $w^{(k)}$ of the equation 
\be\label{uri2} 
(A+(ki)I)w^{(k)}=p^{(k)}. 
\end{equation} 
\end{enumerate} 
 
For any $g\in X$ we have obviously $p^{(k)}\in D(A),\;w^{(k)}\in 
D(A^2)$.
 
In our concrete case the evaluation of (\ref{uri1}) and (\ref{uri2}) is equivalent to the solution of two boundary value problems  
\be\label{irin1} 
\sigma \frac{d^2 p^{(k)}(x)}{dx^2}- 
\mi k  p^{(k)}(x) + c(x) p^{(k)}(x)= g(x), \quad  
p^{(k)}(0)=p^{(k)}(1)=0,\quad k=1, 2, \ldots, 
\end{equation} 
\be\label{irin2} 
\sigma\frac{d^2 w^{(k)}(x)}{dx^2}+\mi k w^{(k)}(x) +c(x)w^{(k)}(x)=  
p^{(k)}(x), \quad 
w^{(k)}(0)=w^{(k)}(1)=0,\quad k=1, 2, \ldots .
\end{equation}

The procedure  stated above for  computing the approximation 
$v_{n, \ell}(t)$ of $v(t)$  makes possible to explore different discretization methods. In particular, we have considered two  
schemes. The first scheme is purely numeric, whereas the second one uses numeric-symbolic computations.
\begin{enumerate}
\item The  first, simplest approach is to discretize the operator $A$ using finite differences methods.  By using finite differences in space with
equispaced points $x_i=ih=i/(N+1)$, $0\leq i\leq N+1$, $h=1/(N+1)$,  for the discretization of the second derivative
we obtain the   discrete analogue $\widehat A$ of the operator $A$, 
\[
\widehat A=\left(\diag[c(x_1) , \ldots, c(x_N)] + \frac{\sigma}{h^2} \left[\begin{array}{cccc}-2 & 1\\ 1& \ddots & \ddots \\& \ddots & \ddots & 1\\ & & 1& -2\end{array}\right] \right) \in \mathbb R^{N\times N}. 
\]
Replacing  $A$ with $\widehat A$ in the procedure for computing $v_{n, \ell}(t)$   is formally equivalent to apply our method for the solution of the first order system  
\begin{equation}\label{mod1}
  \odv{ u}{t}= \widehat A u(t),
  \end{equation}
  with conditions
  \begin{equation}\label{mod1boundary}
  \displaystyle\frac{1}{T}\int_0^Tu(t)dt= f, 
\end{equation}
where  
$ f=\left[ f(x_1), \ldots, f(x_N)\right]^T$.  The  computational problem is of 
the form given in  \eqref{l1},\eqref{l2}   with $\widehat A\in \mathbb R^{N\times N}$.   This finite dimensional first order system  \eqref{mod1}, \eqref{mod1boundary} can also be obtained by  exploiting the classical method of lines \cite{SCH} for solving \eqref{al}, \eqref{lor}.
\item  The second approach combines symbolic and numeric methods.  The evaluation of the derivatives of the function $f$  is performed symbolically.   The computation of \eqref{l1k} is also carried out symbolically. The solution of  (\ref{irin1}) and (\ref{irin2}) can be found by means of  existing software for the solution of BVP's. In our implementation we used the MATLAB function {\tt bvp4c} derived from \cite{Kie} for solving  these  two-point boundary value problems.  The code is  a finite difference algorithm  that implements the three-stage Lobatto IIIa formula. This is a collocation formula and the collocation polynomial provides a $\mathcal C^1$-continuous solution that is fourth order accurate uniformly in $[a,b]$. Mesh selection and error control are adaptively  based on the residual of the continuous solution.
\end{enumerate}

Synthetic  computational problems  have been  designed 
to test  these two approaches.   Our test suite consists of the following two problems: 

\begin{enumerate}
\item The first model example is 
\be\label{alm1}
\left\{\begin{array}{ll}
u_t= \sigma u_{xx} ,&0<x<1,\;0<t< 2\pi,\\
u(0,t)=u(1,t)=0,&0<t< 2\pi.\end{array}
\right.
\end{equation}
\be\label{lorm1}
\frac1{2\pi}\int_0^{2\pi}u(x,t)\;dt=f(x)
\end{equation}
with the data
$$
f(x)=12\frac{1-\displaystyle e^{-2\pi \sigma (3 \pi)^2}}{2 \pi \sigma (3 \pi)^2}\sin(3\pi x) - 7\frac{1-\displaystyle e^{-2\pi \sigma (2 \pi)^2}}{2 \pi \sigma (2 \pi)^2}\sin(2\pi x).
$$
Here the exact solution is 
\[
u(x,t)=12 \displaystyle e^{-\sigma (3 \pi)^2 t}\sin(3\pi x) -7 \displaystyle e^{-\sigma (2 \pi)^2 t}\sin(2\pi x) .
\]

\item The second test is concerned with the differential problem of the form  
\be\label{alm}
\left\{\begin{array}{ll}
u_t=u_{xx} + (4 \pi^2 -1)u(x),&0<x<1,\;0<t< 2\pi,\\
u(0,t)=u(1,t)=0,&0<t< 2\pi.\end{array}
\right.
\end{equation}
\be\label{lorm}
\frac1{2\pi}\int_0^{2\pi}u(x,t)\;dt=f(x)
\end{equation}
with the data
$$
f(x)=\frac{1-e^{-2\pi}}{2 \pi}\sin(2\pi x).
$$
Here the exact solution is 
\[
u(x,t)=\displaystyle e^{-t}\sin(2\pi x).
\]
\end{enumerate}
In the first experiments we tested the accuracy of the first approach based on the semidiscretization in space. Here the error of the computed approximation $v_{n, \widehat \ell}(t)$ is measured with respect to the exact values of the solution function  $u(x,t)$ by 
 \[
 err=\max_{t_i\in \mathcal I_L} \parallel  \widehat u(t_i)-v_{n, \widehat \ell}(t_i)\parallel_\infty/\parallel  \widehat u(t_i)\parallel_\infty , 
 \]
 where $\widehat u(t)=\left[\widehat u_1(t), \ldots, \widehat u_N(t)\right]^T$, $\widehat u_i(t)=u(x_i,t)$, $1\leq i\leq N$.
The coarse and fine grid of points are defined by  $\mathcal I_S=\{t_i=2\pi(i-1)/(m-1)\colon 1\leq i\leq m\}$  and $\mathcal I_L=\{t_i=2\pi(i-1)/(4(m-1))\colon 1\leq i\leq 4(m-1)+1\}$, respectively.
In Table  \ref{table200}  we show the results obtained for the first model problem with $\sigma=1.0e\!-6$, $N=1000$ and $m=100$.  The results are generated by \textbf{Algorithm 1} with $tol=1.0e\!-7$ for different values of $n$.  We report the value of $err$ and 
the minimum and the maximum values of  the length $\ell_i$  of the rational approximation over $t_i\in \mathcal I_s$ as functions of the  degree $n$. 
  \noindent
\begin{table}
		\begin{center}
			\begin{tabular}{|l|l|l|l|} \hline
			
			$n$ & $0$ &$ 1$ & $2$  \\ 	\hline
			$\ell_{min}$ & $7798$ & $177$ & $34$  \\  \hline 
   $\ell_{max}$ & $32716$ & $223$ & $46$  \\ \hline
   $err$ & $1.8e\!-9$ & $1.8e\!-9$ & $1.8e\!-9$  \\ \hline
			\end{tabular}
\end{center}
                \caption{Computed error and the minimum and the maximum values of  the length $\ell_i$  of the rational approximation over $t_i\in \mathcal I_s$ as function of the  degree $n$ for   the first model problem with $\sigma=1.0e\!-6$, $N=1000$ and $m=100$. 
                \label{table200}} 
\end{table}
The discretization error is of order $\sigma h^2 \max|u_{xxxx}|/12$$\simeq 1.0e\!-8$.  The shifted linear systems \eqref{shiftsys} are well conditioned.  Increasing the value of $tol$ does not improve the computed error, which is in accordance with  a priori error bounds. In Table  \ref{table300}  we repeat the same test with $\sigma=0.01$. The discretization error is  now of order $1.0e\!-4$.  The error estimates change significantly since 
the scaling factor $\sigma(N+1)^2$ 
highlights the pathologies  observed in the previous section (compare with the  results  of Test 3).  In particular, the use of  fast expansions with increasing $n$ becomes unreliable.  For $n=2$  \textbf{Algorithm 1} stops but the returned error is larger than $1$.
\noindent
\begin{table}
		\begin{center}
			\begin{tabular}{|l|l|l|l|} \hline
			
			$n$ & $0$ &$ 1$ & $2$  \\ 	\hline
			$\ell_{min}$ & $9093$ & $112$ & $21$  \\  \hline 
   $\ell_{max}$ & $25977$ & $284$ & $33$  \\ \hline
   $err$ & $2.2e\!-4$ & $6.4e\!-6$ & $\star$  \\ \hline
			\end{tabular}
\end{center}
                \caption{Computed error and the minimum and the maximum values of  the length $\ell_i$  of the rational approximation over $t_i\in \mathcal I_s$ as function of the  degree $n$ for   the first model problem with $\sigma=1.0e\!-2$, $N=1000$ and $m=100$.  For $n=2$ the computed error is larger than 1. 
                \label{table300}} 
\end{table}

\begin{figure}
  \centering
  \subfloat[]{\includegraphics[width=0.5\textwidth]{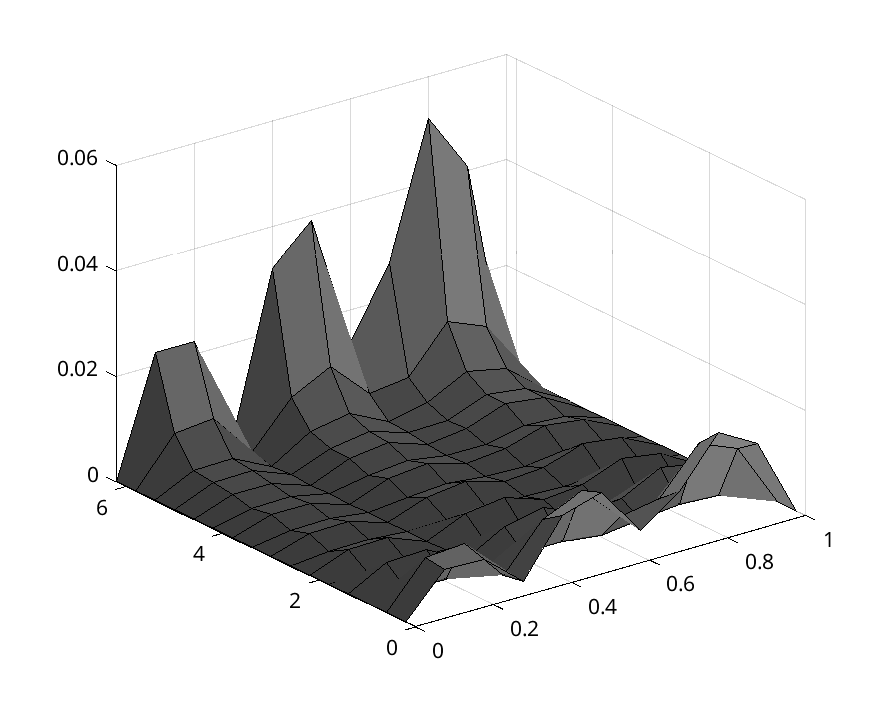}}
  \hfill
  \subfloat[]{\includegraphics[width=0.5\textwidth]{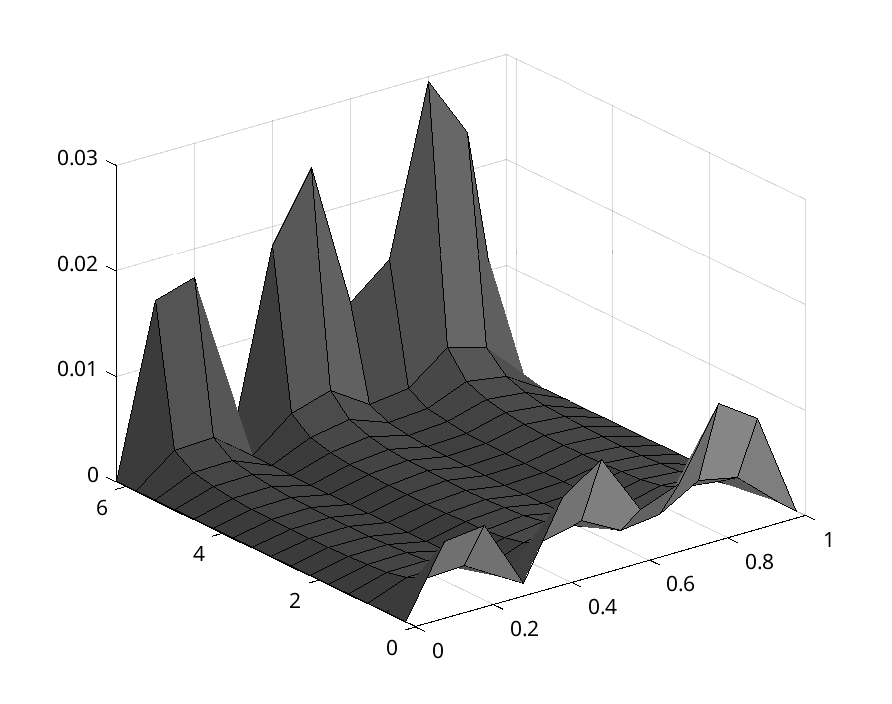}}
  \caption{Plot of $err(x,t)$ generated by solving the first model problem  with $\sigma=0.01$ for $n=0$ with $\ell=200$  in $(a)$ and $\ell=400$ in $(b)$.}
\label{ff500}
\end{figure}

\begin{figure}
  \centering
  \subfloat[]{\includegraphics[width=0.5\textwidth]{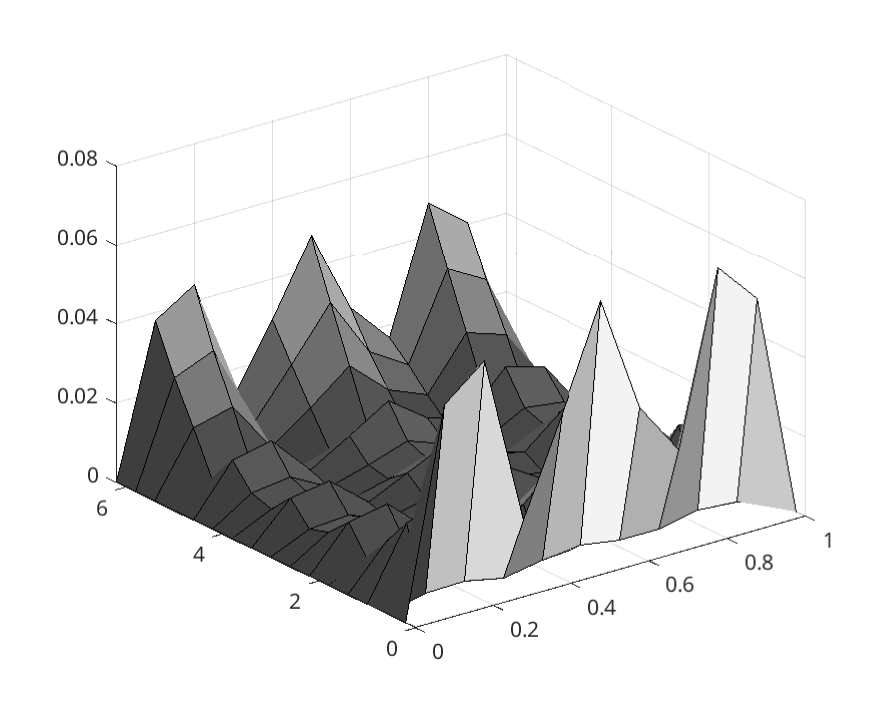}}
  \hfill
  \subfloat[]{\includegraphics[width=0.5\textwidth]{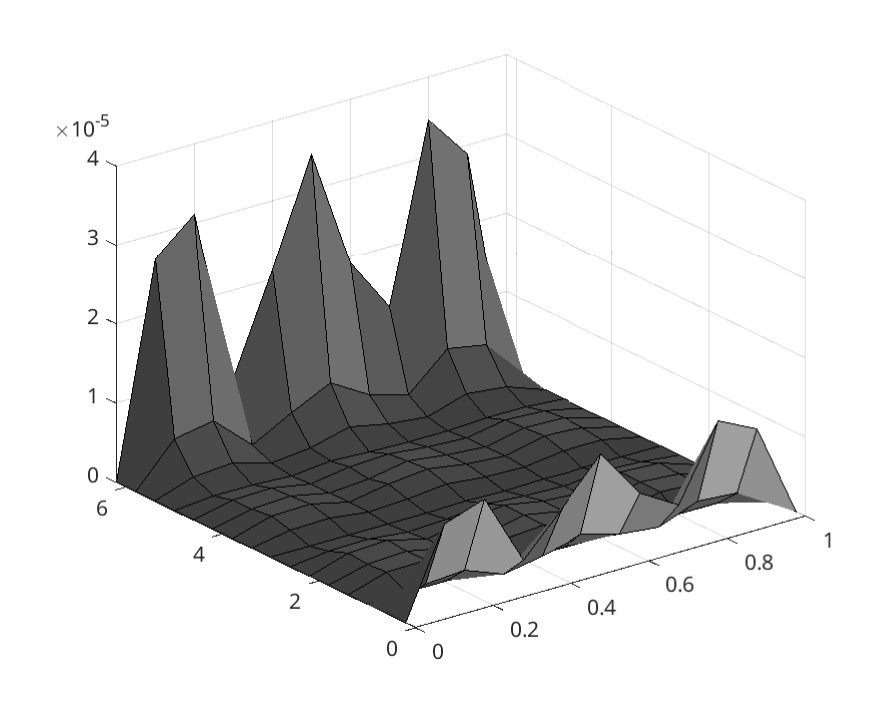}}
  \caption{Plot of $err(x,t)$  generated by solving the first model problem  with $\sigma=0.01$ for $n=1$ with $\ell=10$  in $(a)$ and $\ell=100$ in $(b)$. }
\label{ff600}
\end{figure}

\begin{figure}
  \centering
  \subfloat[]{\includegraphics[width=0.5\textwidth]{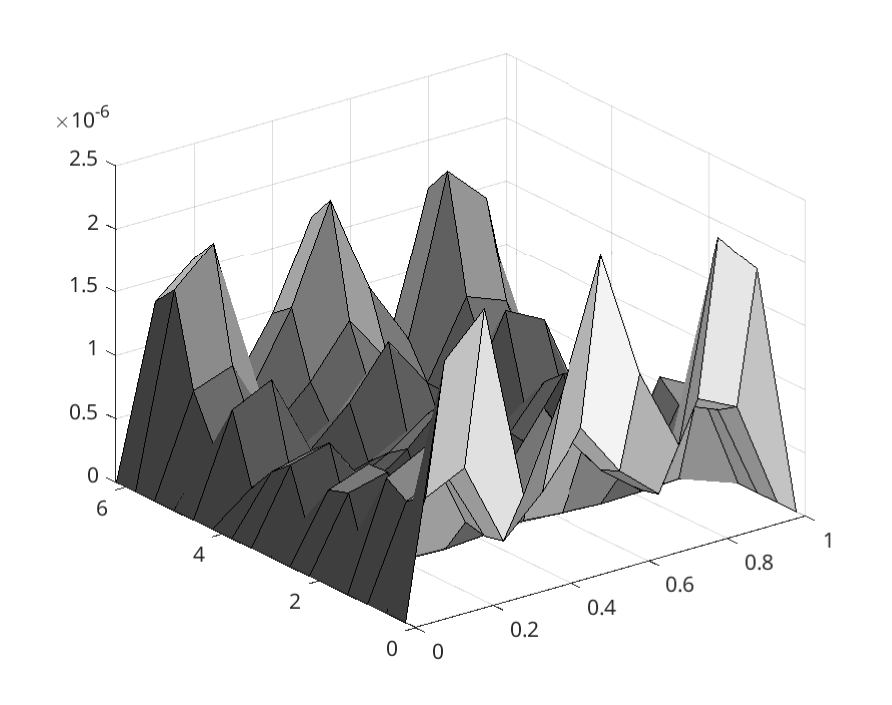}}
  \hfill
  \subfloat[]{\includegraphics[width=0.5\textwidth]{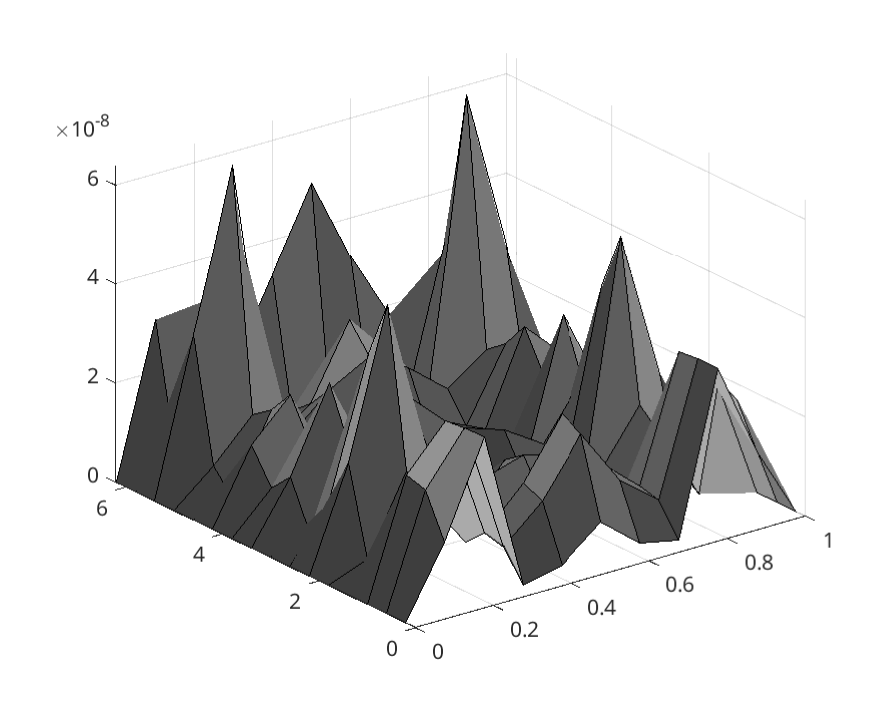}}
  \caption{Plot of $err(x,t)$ generated by solving the first model problem  with $\sigma=0.01$ for $n=2$ with $\ell=10$  in $(a)$ and $\ell=20$ in $(b)$.}
\label{ff700}
\end{figure}

\begin{figure}
  \centering
  \subfloat[]{\includegraphics[width=0.5\textwidth]{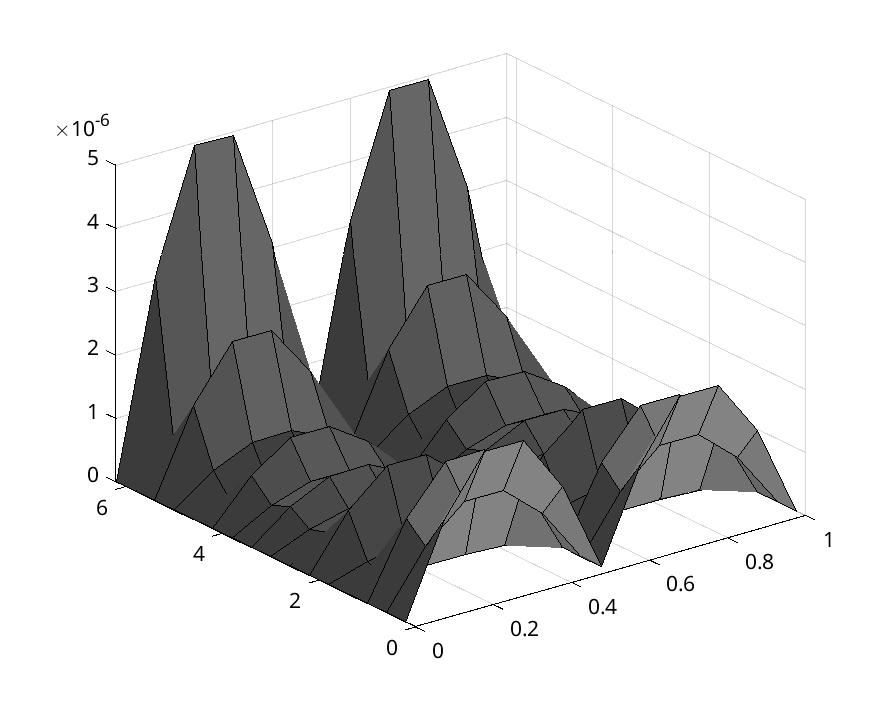}}
  \hfill
  \subfloat[]{\includegraphics[width=0.5\textwidth]{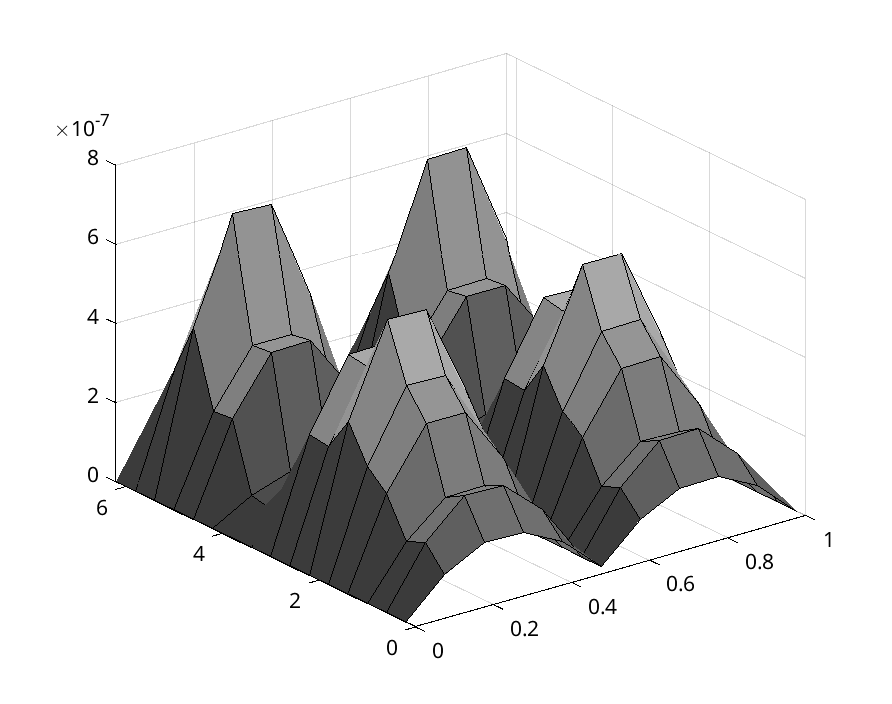}}
  \caption{Plot of $err(x,t)$ generated by solving the second  model problem   for $n=1$ with $\ell=20$  in $(a)$ and $\ell=40$ in $(b)$.}
\label{ff800}
\end{figure}

The second approach is based on computing the expansion of $v(t)$   as given in Theorem \ref{main} for the infinite dimensional operator $A=\sigma \frac{d^2}{dx^2}+c(x)$. Recall that in our implementation we used the MATLAB function {\tt bvp4c} for solving  the  two-point boundary value problems \eqref{irin1}, \eqref{irin2}.
The number of boundary value problems to be solved depends on the value of $\ell$ in \eqref{sch1} which is fixed a priori. The output is an interpolating function $v_{n, \ell}(x,t)$.  For the sake  of graphical illustration the  computed error  is  shown by plotting the function $err(x,t)=|v_{n, \ell}(x,t)-u(x,t)|$ over  the 
domain $[0,1]\times [0,2 \pi]$.  In Figure \ref{ff500}, \ref{ff600}
and \ref{ff700} we show the plots  generated by solving the first model problem with $\sigma=0.01$ for different values of $n$ and $\ell$. 
In Figure \ref{ff800} we describe our results for the second model problem with $n=1$ and $\ell=20$ and $\ell=40$. 

In summary, our preliminary experience with the infinite dimensional (functional) approach indicates that, when applicable, expansions with increasing $n$  allow the computation of quite accurate  approximations of $u(x,t)$ already for small or moderate values of the truncation level $\ell$.

\section{Conclusion and Future Work}\label{sec:conclusions}
In this paper we have devised a  mixed polynomial/rational expansion for the  solution of an abstract first-order differential problem with nonlocal conditions. We have proposed an algorithm that exploits this expansion to approximate numerically the solution and tested it on several examples, both in finite- and infinite-dimensional settings. Of particular interest are implementations where the expansion is directly applied to an infinite-dimensional operator and discretization only intervenes in the solution of the resulting family of shifted boundary value problems.

There are  still many points for future developments. Specifically:
\begin{enumerate}
\item The derivation of precise residual estimates for this expansion also related with the design of efficient stopping criteria. 
\item  The study  of alternative acceleration schemes not involving  high-degree Bernoulli polynomials  in order to to relax the regularity assumptions.
\item The design of  efficient numerical methods for solving the sequence of shifted boundary value problems \eqref{irin1}, \eqref{irin2}. 
\item The application of our proposed method to more general differential problems such as the two-dimensional parabolic equation and the biharmonic equation. 
\end{enumerate}

We also point out that, for suitable classes of problems, one may apply the mixed polynomial/rational expansion presented here to devise hybrid approaches that combine symbolic/functional and purely numerical techniques. As an example, consider Problem 1 in Section 3, i.e., \eqref{alm1}, \eqref{lorm1}. Starting from the finite-dimensional approach based on semi-discretization in space, we may incorporate more and more symbolic elements to circumvent potential sources of instability and improve the quality of the results. For instance, since the function $f(x)$ for \eqref{alm1}, \eqref{lorm1} is known analytically and the action of matrix $A$ corresponds to a second derivative, one may pre-compute symbolically $f''(x)$ and $f^{(4)}(x)$, evaluate them on the points of the space grid, and use the resulting vectors in place of $Af$ and $A^2f$ in the polynomial/rational expansion. Moreover, we may also try evaluating the polynomial part of the expansion symbolically, if the polynomial degree is small (but high enough to be numerically troublesome). Preliminary experiments show that, with these modifications, the test presented in Table \ref{table300} for $n=2$ yields an error that is no longer larger than 1, but comparable to case $n=1$, with a smaller number of rational terms.

\bibliographystyle{plain}
\bibliography{newton}

\end{document}